\documentclass{article}

\usepackage{graphicx}%
\usepackage{multirow}%
\usepackage{amsmath,amssymb,amsfonts}%
\usepackage{mathrsfs}%
\usepackage[title]{appendix}%
\usepackage{xcolor}%
\usepackage{textcomp}%
\usepackage{manyfoot}%
\usepackage{booktabs}%
\usepackage{algorithm}%
\usepackage{algorithmicx}%
\usepackage{algpseudocode}%
\usepackage{listings}%

\newtheorem{remark}{Remark}

\usepackage[style=numeric]{biblatex}
\title{A simple algorithm for the summation of alternating series}

\author{Eleonora Denich\footnotemark[2]
        \and Paolo Novati \footnotemark[3] \and Alvise Sommariva \footnotemark[4]}

\begin{document}

\maketitle

\renewcommand{\thefootnote}{\fnsymbol{footnote}}

\footnotetext[2]{Department of Mathematics, Informatics and Geosciences, University of Trieste, via Valerio 12/1, 34127, Trieste, Italy ({\tt eleonora.denich@units.it}).}
\footnotetext[3]{Department of Mathematics, Informatics and Geosciences, University of Trieste, via Valerio 12/1, 34127, Trieste, Italy ({\tt novati@units.it}).}
\footnotetext[4]{Department of Mathematics "Tullio Levi-Civita", University of Padova, Via Trieste, 63, 35121 Padova, Italy ({\tt alvise@math.unipd.it}).}

\begin{abstract}
This paper deals with the computation of the sum of alternating series, whose general terms can be expressed by means of analytic functions. After rewriting the series as a weighted integral with the Abel weight, we employ the sinc (trapezoidal) rule and analyze the remainder term with respect to the number of quadrature points. We provide some numerical experiments to show the reliability of the derived error estimate and to test the algorithm for automatic summation with prescribed accuracy.
All the Matlab codes used in the present paper can be found as open-source software at \cite{GITDNS}.

%\textcolor{black}{ALVISE: Nell'articolo parliamo sempre dei trapezi ma poi lo chiamiamo {\it{sinc}}. Che sia da dire qualcosa?}

\end{abstract}


{\bf{AMS}}:
Primary 40A25; Secondary 65D30, 65D32, 65G99.

\section{Introduction}

Given a function $f$, many numerical methods for the fast summation of series of the type
\begin{equation} \label{sum1}
    \sum_{j=0}^{+\infty} f(j)
\end{equation}
or
\begin{equation} \label{sum2}
    \sum_{j=0}^{+\infty} (-1)^j f(j)
\end{equation}
are based on the reformulation of the problem in terms of the computation of weighted integrals (clearly, $f$ has to satisfy suitable hypotheses, that depend on the method).
Examples are given by the famous Abel \cite{Abel_2012}, Plana \cite{PLANA} and Lindelhöf \cite{L1905} formulas, but other interesting approaches have been introduced in the past, see e.g. \cite{GM85} and the contribute given by Dahlquist in \cite{DI97,DII97,DIII99} for a wide overview.
The methods based on the computation of weighted integrals in general require Gaussian quadrature on the positive real line with a (non classical) weight function that decays exponentially.
In particular, in \cite{GM85} quadrature rules on $(0,+\infty)$ for the Einstein and Fermi weight functions
\begin{equation*}
    \epsilon(t) = \frac{t}{e^t-1}, \quad \varphi(t) = \frac{1}{e^t+1}
\end{equation*}
have been constructed and applied respectively to series of type (\ref{sum1}) and (\ref{sum2}), whose general terms are expressible as derivatives of a Laplace transform.
In \cite{M94}, the author gives an alternative approach based on the contour integration over a rectangle in the complex plane and on the reduction of (\ref{sum1}) and (\ref{sum2}) to a quadrature problem on $(0,+\infty)$, with respect to the weight functions
\begin{equation*}
    w_1(t) = \frac{1}{\cosh^2t} \quad {\rm and} \quad w_2(t) = \frac{\sinh t}{\cosh^2t},
\end{equation*}
respectively.
More recently, in \cite{DJ25}, the author studies orthogonal polynomials with respect to the Abel weight function 
\begin{equation*}
    \omega(t) = \frac{t}{2 \sinh(\pi t)}. 
\end{equation*} 
Next, by exploiting the Abel summation formula \cite{Abel_2012}, the corresponding Gaussian rule has then been used for summing the infinite alternating series (\ref{sum2}), where $f$ is assumed to be holomorphic in an open region containing the half plane $\Re(z) \geq 0$.
All these methods are generally quite effective but, sometimes, non trivial additional information on $f$ are required, as for instance its primitive (\cite{M94}) and its inverse Laplace transform (\cite{GM85}).
Moreover, up to our knowledge, for none of them the error analysis is available.

In this work we concentrate on alternating series of type (\ref{sum2}) and analyze the sinc rule applied to Abel's formula.
Moreover, we present a quite accurate error analysis that allows to construct an automatic procedure for the computation of the alternating series with a prescribed accuracy.

This work is organized as follows.
In Section \ref{section 2} we recall the basic results concerning the trapezoidal (sinc) rule for evaluating integrals over $(-\infty,+\infty)$.
In Section \ref{section 3} we present the method used for evaluating (\ref{sum2}), together with the error analysis. Some numerical experiments are given in Section \ref{section 4}.
In the appendix, we describe the Matlab codes used in this work, available as open-source software at {\cite{GITDNS}}.

\section{General results} \label{section 2}

In this section we recall some theoretical results concerning the trapezoidal approximation
\begin{equation*}
    I(G) = \int_{-\infty}^{+\infty} G(x) dx \approx h \sum_{k=-\infty}^{+\infty } G(kh),
\end{equation*}
where $G \colon \mathbb{R} \rightarrow \mathbb{R}$ and $h>0$. Denoting by
\begin{equation*}
    T_{N,M,h} (G) = h \sum_{k=-N}^M G(kh),
\end{equation*}
the truncated trapezoidal rule, with $N,M$ positive integers, we have that for the quadrature error 
\begin{equation*}
    \mathcal{E}_{N,M,h} = \left\vert I(G)-T_{N,M,h}(G) \right\vert
\end{equation*}
it holds
\begin{equation*}
    \mathcal{E}_{N,M,h} \leq \mathcal{E}_D+\mathcal{E}_{T_L}+\mathcal{E}_{T_R},
\end{equation*}
where 
\begin{equation} \label{E_D}
    \mathcal{E}_{D}= \left\vert \int_{-\infty}^{+\infty} G(x) dx-h \sum_{k=-\infty}^{+\infty}G(kh) \right\vert
\end{equation}
is the discretization error and 
\begin{equation} \label{E_T}
    \mathcal{E}_{T_L}= h \left\vert \sum_{k=-\infty}^{-N-1}G(kh) \right\vert, \quad \mathcal{E}_{T_R}= h \left\vert \sum_{k=M+1}^{+\infty}G(kh) \right\vert
\end{equation}
are the truncation errors.

We assume that the function $G$ restricted to the real axis is real valued and that it is analytic in an infinite strip domain containing the real axis, except for a pair of simple isolated singularities $z_0$ and its conjugate $\overline{z_0}$.
In this setting, in \cite{DN24} it has been shown that
\begin{equation}\label{50}
    \mathcal{E}_D \lesssim 4 \pi \vert \rho_0 \vert e^{-2d \frac{\pi}{h}},
\end{equation}
where $d = \Im(z_0)$ and $\rho_0 = {\rm Res} (G,z_0)$ is the residue of $G$ at $z_0$ (we have implicitly assumed that $\Im(z_0)>0$).
The symbol $\lesssim$ indicates less than or asymptotically equal.
Then, in order to optimize the rate of convergence with respect to the total number of quadrature points $M+N+1$, one has to impose that $\mathcal{E}_D, \mathcal{E}_{T_L},\mathcal{E}_{T_R}$ have the same exponential decay. 
In this way, one determines the values of $h,N,M$ and use them to obtain an approximation of the total error $\mathcal{E}_{N,M,h}$.

We observe that, if $G$ is an even function, then $I(G) \approx 2 T_{M,h}(G)$, where
\begin{equation} \label{formula C}
    T_{M,h}(G) = h \left(\frac{1}{2}G(0)+\sum_{k=1}^M G(kh)\right).
\end{equation}
In this situation, {\textcolor{black}{it can be easily checked that}} 
\begin{equation*}
    \mathcal{E}_{M,h} = \left\vert \int_0^{+\infty}G(x) dx -T_{M,h}(G) \right\vert \leq \frac{\mathcal{E}_D}{2}+\mathcal{E}_{T_R},
\end{equation*}
where $\mathcal{E}_D$ and $\mathcal{E}_{T_R}$ are defined in (\ref{E_D}) and (\ref{E_T}).

%{\textcolor{black}{ALVISE (10.8.26): questa parte mi e' un po' opaca; secondo me $T_{M,M,h}(G):=T_{M,h}(G)=h (G(0)+2\sum_{k=1}^M G(kh))$. Non capisco al momento perch\'e si dimezzino gli errori, devo pensarci un po' su, guardo il vostro articolo. Mi verrebbe da dire che l'errore \'e  $\mathcal{E}_{M,h} \leq {\mathcal{E}_D}+2\mathcal{E}_{T_R}.$}}

%\textcolor{blue}{ELEONORA 17.08: ho corretto la formula per $T_{M,h}$ (in effetti mancava l'$h$ davanti). C'era un altro errore prima della formula di $T_{M,h}$, mancava un 2 (l'ho corretto in blu). Ora dovrebbe più chiaro il perché gli errori si dimezzano. Comunque, abbiamo rifatto i conti e la maggiorazione di $\mathcal{E}_{M,h}$ è corretta.}

\section{The trapezoidal rule for the computation of (\ref{sum2})} \label{section 3}

Let $f$ be a function analytic in a region of the complex plane containing the half plane $\Re (z) \geq 0$ and such that (see \cite[p.18]{M2017})
\begin{enumerate}
    \item $\lim_{\vert y \vert \rightarrow +\infty }e^{-\vert 2 \pi y \vert}\vert f(x \pm iy)\vert=0$ uniformly in $x$ on every finite interval,
    \item $\int_{-\infty}^{+\infty}\vert f(x+iy)-f(x-iy) \vert e^{-\vert 2 \pi y \vert} dy$ exists for every $x \geq 0$ and tends to zero when $x \rightarrow + \infty$.
\end{enumerate}
%{\textcolor{black}{ALVISE 12.08.26: a me sembra che la referenza giusta coi trapezi sia {\cite[p.67]{DIII99}} per $m=0$.}}
%\textcolor{blue}{ELEONORA 17.08: ho corretto la referenza prima della due condizioni sulla $f$ e quella per la formula di Abel.}
Then, the Abel formula (see also \cite[p.67]{DIII99}) gives
\begin{equation*}
    A(f) := \sum_{j=0}^{\infty} (-1)^j f(j) = \frac{1}{2}f(0)-\int_0^{+\infty} \frac{f(ix)-f(-ix)}{2i \sinh (\pi x)}dx. 
\end{equation*}
Assuming $\overline{f(z)}=f(\overline{z})$, that is, the restriction of $f$ to the real numbers is real valued, and since $x \in [0, +\infty)$, the above relation can be rewritten as
\begin{equation} \label{formula D}
    A(f) = \frac{1}{2}f(0)-\int_0^{+\infty}\frac{\Im(f(ix))}{\sinh (\pi x)}dx.
\end{equation}
%{\textcolor{black}{ALVISE 11-08-26. Che sia da dire che la richiesta sui coniugi e' automaticamente verificamente per funzioni real-analytic?}}
%\textcolor{blue}{ELEONORA 17.08: sembra essere una cosa nota ma se vuoi possiamo aggiungere.}
%\textcolor{blue}{ALVISE 18.08: direi di scriverlo, il paper potrebbe finire in mano a qualcuno che sa di quadratura ma non cos\'i bene di analisi complessa.}
Denoting by 
\begin{equation*}
    s_f(x) = \frac{\Im(f(ix))}{\sinh (\pi x)}, \quad x \in \mathbb{R},
\end{equation*}
and observing that $s_f(x) = s_f(-x)$, the trapezoidal rule for evaluating $A(f)$ leads to
\begin{equation*}
    A(f)=\frac{1}{2}f(0)-h \left( \frac{s_f(0)}{2}+\sum_{k=1}^M s_f(kh) \right) - \mathcal{R}_M(s_f),
\end{equation*}
where 
\begin{equation*}
    \mathcal{R}_M(s_f)= \int_0^{+\infty}s_f(x) dx -T_{M,h}(s_f)
\end{equation*}
(see (\ref{formula C})).
We remark that the function $s_f$ has a removable singularity in $0$ since
\begin{equation} \label{eq1}
    \frac{\Im(f(ix))}{x} \rightarrow f'(0),
\end{equation}
for $x \rightarrow 0$ (complex step derivative approximation, see \cite{MSA}).
%{\textcolor{black}{ALVISE 11-08-26. Che sia da dire che la parte immaginaria di una funzione analitica \'e analitica (alla fine lo si manda a una rivista numerica e qualcuno potrebbe non essere esperto)?}}
%\textcolor{blue}{ELEONORA 17.08: non credo sia vero. Comunque, non è necessario perché si lavora con una funzione a valori reali.}
Therefore, since $\lim_{x \rightarrow 0} {\sinh(x)}/{x}=1$, we can define
\begin{equation*}
    s_f(0) = \frac{f'(0)}{\pi}.
\end{equation*}
An interesting feature is that the approximation $f'(0) \approx {\Im(f(i\delta))}/{\delta} $ suggested by (\ref{eq1}) is stable even with $\delta \approx 1e-100$, a property that is particularly appealing when the derivative of $f$ in $0$ is not available.
Thus,
\begin{equation}\label{83}
    A(f)=\frac{1}{2}f(0)-h \left( \frac{\Im(f(i\delta))}{2 \delta \pi}+\sum_{k=1}^M s_f(kh) \right) - \mathcal{R}_M(s_f),
\end{equation}
where have considered negligible the error of the complex step derivative approximation.

Before going on with our analysis, we make some considerations.

First, in view of (\ref{E_D}) and (\ref{50}), for $h$ fixed the rate of convergence is mainly determined by the magnitude of $d$. 

Next, since $s_f(x) = {\Im(f(ix))}/{\sinh (\pi x)}$, the function $\sinh(\pi z)$ introduces the poles
\begin{equation*}
    z_{\ell}=i \ell, \quad \ell \in \mathbb{Z}.
\end{equation*}
Being $z_0=0$ a removable singularity, the poles closest to the real axis are $z_1=i$ and its conjugate $z_{-1}=-i$.
On the other side, each singular point of $f$, say $w=a+ib$, $a<0$, determines a couple of singularities $ w_{\pm}=\pm b \mp ia$, for $f(iz)-f(-iz)$.
As a consequence, $d = \min\lbrace \vert a \vert,1\rbrace$.
Clearly, the best situation is when $d=1$ and this is always possible by discarding the first term of the sum.

Indeed, setting ${\tilde{f}}(x)=f(x+1)$ and
\begin{equation} \label{s tilde}
    s_{\tilde{f}}(x) = \frac{\Im({\tilde{f}}(ix))}{\sinh (\pi x)}=\frac{\Im(f(ix+1))}{\sinh (\pi x)}, \quad x \in \mathbb{R},
\end{equation}
from
\begin{eqnarray} 
A({\tilde{f}}) &=& \frac{1}{2} {\tilde{f}}(0)-h \left( \frac{\Im({\tilde{f}} (i\delta))}{2 \delta \pi}+\sum_{k=1}^M {s_{\tilde{f}}}(kh) \right) - \mathcal{R}_M({s_{\tilde{f}}}) \nonumber \\
&=& \frac{1}{2} {{f}}(1)-h \left( \frac{\Im({{f}} (i\delta+1))}{2 \delta \pi}+\sum_{k=1}^M {s_{\tilde{f}}}(kh) \right) - \mathcal{R}_M({s_{\tilde{f}}})
\end{eqnarray}
we get
\begin{eqnarray}\label{formulaA(f)}
    A(f)&=&\sum_{j=0}^{\infty} (-1)^j f(j) = f(0) - \sum_{j=0}^{\infty} (-1)^jf(j+1) = f(0)-A({\tilde{f}}) \nonumber \\
    &=& f(0)-\left[ \frac{1}{2}f(1)-h \left( \frac{\Im(f(i\delta+1))}{2 \delta \pi}+\sum_{k=1}^M {s_{\tilde{f}}}(kh) \right) - \mathcal{R}_M({s_{\tilde{f}}}) \right].
\end{eqnarray}
Consequently, we can consider the formula
\begin{equation} \label{metodo 2}
    A(f) = f(0)-\frac{1}{2}f(1) +h \left( \frac{\Im(f(i\delta+1)}{2 \delta \pi}+\sum_{k=1}^M s_{\tilde{f}}(kh) \right) + \mathcal{R}_M({s_{\tilde{f}}}).
\end{equation}
As an effect, since the singularities given by the factor $f(ix+1)-f(-ix+1)$ are $z =\pm b\mp i(a-1)$ with, $a<0$, we obtain $d=1$ and we can use the analysis of Section \ref{section 2} to study the error of method (\ref{metodo 2}).

In particular, we have that
\begin{equation*}
    \mathcal{E}_D \lesssim 2 \pi \left\vert {\rm Res}(s_{\tilde{f}}(z),z_0) \right\vert e^{-2d\frac{\pi}{h}}, \quad h \rightarrow 0,
\end{equation*}
with $z_0 = i$, $d=1$. Since 
$$
{\rm Res}\left(s_{\tilde{f}}(z),z_0\right) =-\frac{f(0)-f(2)}{2 \pi i},
$$ 
one finds
\begin{equation*}
    \mathcal{E}_D \lesssim  \left\vert {f(0)-f(2)} \right\vert e^{-2\frac{\pi}{h}}, \quad h \rightarrow 0.
\end{equation*}
By using 
$$
\frac{1}{\sinh (\pi x) } \sim 2 e^{-\pi x}, \quad x \rightarrow + \infty
$$ 
(the symbol $\sim$ denotes asymptotic equality), assuming that 
\begin{equation} \label{formula A}
    \left\vert \Im (f(i(Mh+1)+1)) \right\vert \leq C,
\end{equation} 
definitely,  for the truncation error $\mathcal{E}_{T_R}$ we obtain
\begin{align}
    \mathcal{E}_{T_R} &= h \left\vert \sum_{k=M+1}^{+\infty} s_{\tilde{f}}(kh) \right\vert \lesssim 2Ch \sum_{k=M+1}^{+\infty} e^{-\pi kh}  \\
     &\leq 2 C \int_{Mh}^{+\infty}e^{-\pi t}dt = \frac{2C}{\pi} e^{-\pi M h}. \notag
\end{align}
Now, in order to have the same asymptotic behavior for  $\mathcal{E}_{D}$ and $\mathcal{E}_{T_R}$, we impose
\begin{equation*}
    2 \frac{\pi}{h} = \pi M h,
\end{equation*}
that is,
\begin{equation} \label{formula B}
    h = \sqrt{\frac{2}{M}},
\end{equation}
obtaining
\begin{equation} \label{stima errore}
    \mathcal{E}_{M,h} \leq \frac{\mathcal{E}_D}{2}+ \mathcal{E}_{T_R} \lesssim \left( \vert f(0)-f(2) \vert +\frac{2C}{\pi} \right)e^{-\pi \sqrt{2M}}, \quad M \rightarrow +\infty.
\end{equation}

\begin{remark}
    We remark that in formula (\ref{3.2 bis}) we have not considered the decay rate of $\Im(f(ix+1))$ (see (\ref{s tilde})), and hence formula (\ref{stima errore}) may be a bit conservative whenever $f(j)$ rapidly goes to $0$ for $j\rightarrow+\infty$.
\end{remark}

%\textcolor{black}{ALVISE 12.08.26: Ho un dubbio. Tutta la teoria sugli errori \'e per integrali sulla retta reale. Ci sono ipotesi. Diciamo qualcosa sulle funzioni pari, ma sempre sulla retta reale. Qui invece facciamo le cose sulla semiretta reale. Come si passa da uno all'altro? Un prolungamento alla retta reale per simmetria? Mantiene tutte le ipotesi? Noto poi una cosa di potenziale interesse, Dahlquist usa pure le midpoint composte.}
%\textcolor{blue}{ELEONORA 17.08: per poter usare la teoria generale $(-\infty, +\infty)$ al nostro caso, è sufficiente estendere la funzione integranda sulla semiretta negativa in maniera che risulti pari.}

\section{Numerical experiments} \label{section 4}

In this section we present some numerical tests in which we show the reliability of the method proposed and error estimate (\ref{stima errore}). When the primitive $F$ of the function $f$ is known, we also compare the trapezoidal rule with the method introduced by Milovanović in \cite{M94}. 
We remind the reader that in the latter, observing that $$A(f)=\sum_{k=0}^{\infty} (-1)^k f(k)=\sum_{k=0}^{l-1} (-1)^k f(k)+\sum_{k=l}^{\infty} (-1)^k f(k),$$
under suitable assumptions on $f$, the author first computes $\sum_{k=0}^{l-1} (-1)^k f(k)$ and then approximates 
\begin{equation*}
    A_l(f)=\sum_{k=l}^{\infty} (-1)^k f(k)
\end{equation*}
by $A_{l,n}$, an evaluation of  a certain integral involving $F$, via a suitable $n$-point Gaussian rule (see the appendix of the present paper for further details).

As numerical tests, we consider the sums 
\begin{equation*}
    A(f_i)=\sum_{k=0}^{+\infty} (-1)^k f_i(k),
\end{equation*}
where
\begin{enumerate}
    \item $f_1(x) = \frac{1}{(x+2)^2+1}$, with primitive $F_1(x)= \arctan(x+2)$;
    \item $f_2(x) = \frac{1}{(x+2)^{\alpha}}$, $\alpha >0$, with primitive $F_2(x) = \frac{(x+2)^{-\alpha+1}}{1-\alpha}$;
    \item $f_3(x) = \frac{1}{(x+1)^{\beta}+1}$, $\beta >0$.
\end{enumerate}
{\textcolor{black}{Though our Matlab codes perform tests with general parameters $\alpha > 0$, $\beta >0$, for illustration purposes we have considered the cases in which $\alpha=1.5$, $\beta=3/2$.}}

As for the reference results of $A(f_i)$, we have used (\ref{formula D}) in which we have numerically computed $\int_0^{+\infty} \frac{\Im(f(ix+1))}{\sinh(\pi x)} dx$ by Matlab {\tt{integral}} routine, with a relative tolerance equal to machine precision, obtaining
\begin{enumerate}
    \item $S^{(1)} \approx 0.136014527491066572$;
    \item $S^{(2)} \approx 0.234852975374592098$;
    \item $S^{(3)} \approx  0.335654070827147765$.
\end{enumerate}

%{\textcolor{black}{ALVISE 12-08-26. Le funzioni $f_j$ verificano le ipotesi richieste da Abel?}}

%\textcolor{blue}{ELEONORA 17.08: non abbiamo fatto i test ma speriamo di si.}

In all the Figures below, we plot the absolute summation errors obtained by
\begin{itemize}
\item the direct summation {\tt{sum}} using $m$ functional evaluations,
\item the trapezoidal approach {\tt{trapz}},
\end{itemize}
increasing the number of functional evaluations of $f_i$ as well as the numerical estimate (\ref{stima errore}), where we have set
    \begin{equation} \label{costante C}
        C = \left\vert \Im(f(i\sqrt{2}+1) \right\vert,
    \end{equation}
see (\ref{formula A}) and relation (\ref{formula B}).

As for the evaluation of $A(f_1)$ and  $A(f_2)$, being a primitive $F_1$ and $F_2$ at hand, we also consider the summation methods proposed by Milovanovic, in which {\tt{mil1}}, {\tt{mil3}}, {\tt{mil5}} represent respectively the quantities $f(0)+A_{1,n}(f)$, $\sum_{k=0}^2 (-1)^k f(k)+A_{3,n}(f)$, $\sum_{k=0}^4 (-1)^k f(k)+A_{5,n}(f)$. We observe that {\tt{mil1}}, {\tt{mil3}}, {\tt{mil5}} require, depending on the method, $n+1$, $n+3$ and $n+5$ functional evaluations.

\begin{figure}
\begin{center}
\includegraphics[width=10cm]{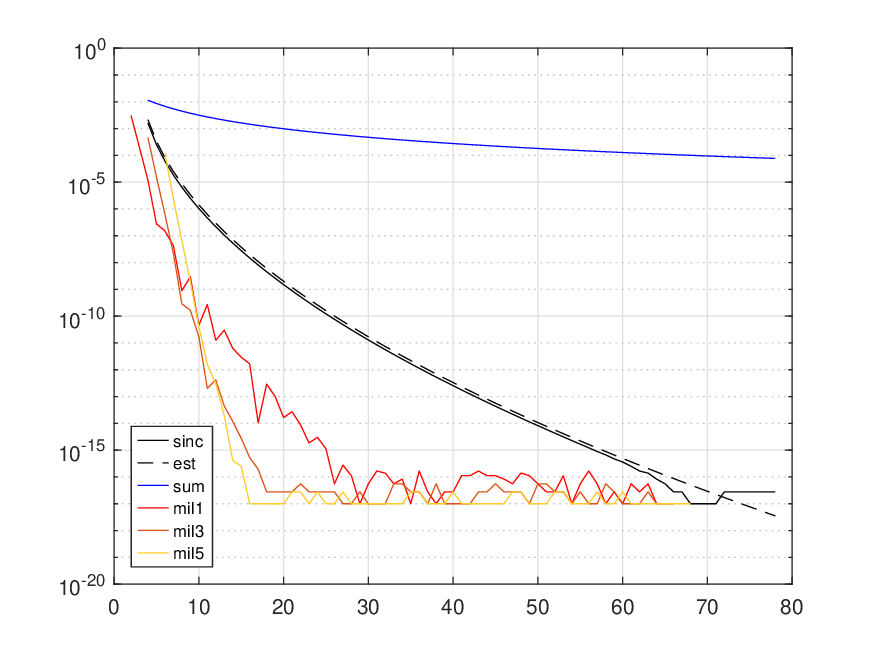}
\caption{Absolute errors with respect to the number of evaluations of the $f_1$ by the direct summation {\tt{sum}}, the trapezoidal rule {\tt{sinc}} (together with the error estimate {\tt{est}}) and the methods {\tt{mil1}}, {\tt{mil3}}, {\tt{mil5}} proposed in \cite{M94}.}
\label{fig1}
\end{center}
\end{figure}

\begin{figure}
\begin{center}
\includegraphics[width=10cm]{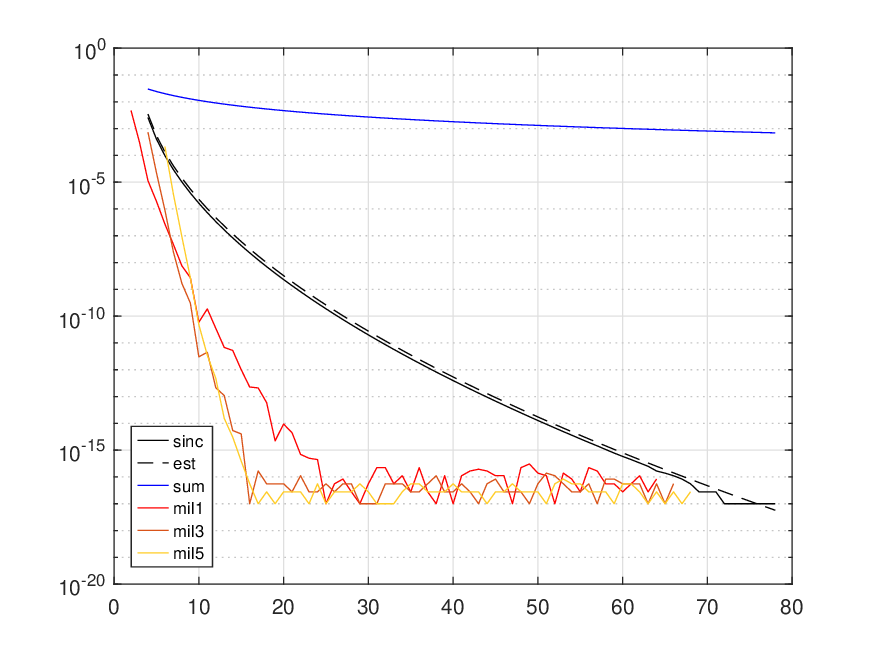}
\caption{Absolute errors with respect to the number of evaluations of the $f_2$ by the direct summation {\tt{sum}}, the trapezoidal rule {\tt{sinc}} (together with the error estimate {\tt{est}}) and the methods {\tt{mil1}}, {\tt{mil3}}, {\tt{mil5}} proposed in \cite{M94}.}
\label{fig2}
\end{center}
\end{figure}

\begin{figure}
\begin{center}
\includegraphics[width=10cm]{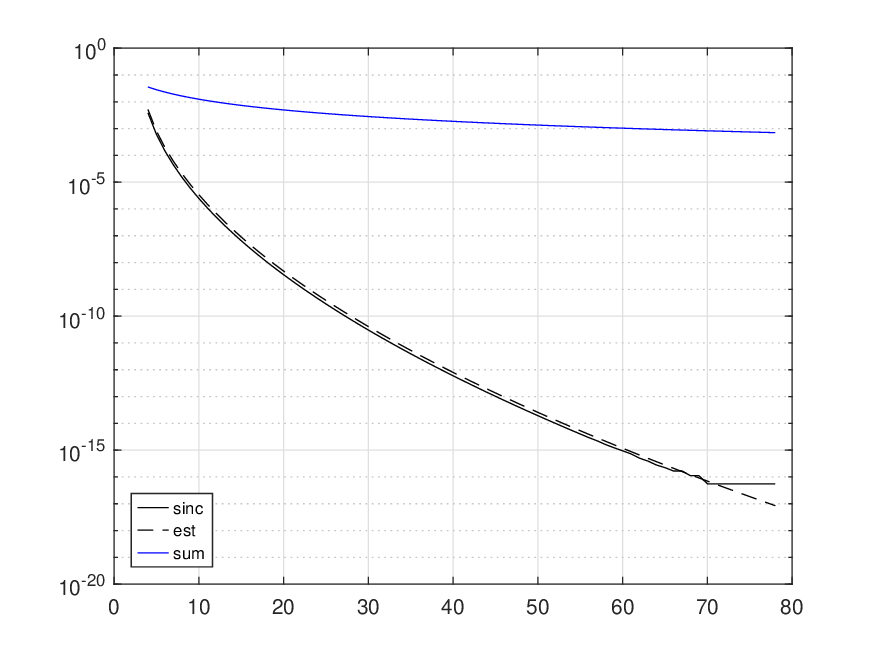} 
\caption{Absolute errors with respect to the number of evaluations of the $f_3$ by the direct summation {\tt{sum}}, the trapezoidal rule {\tt{sinc}} (together with the error estimate {\tt{est}}) and the methods {\tt{mil1}}, {\tt{mil3}}, {\tt{mil5}} proposed in \cite{M94}.}
\label{fig3}
\end{center}
\end{figure}

In all the figures, the series are slowly convergent. In these tests the method based on the trapezoidal rule reaches an absolute error of the order of machine precision, with about $60$ evaluations, and the numerical estimate of the trapezoidal rule error is particularly tight. Milovanovic summation algorithms {\tt{mil1}}, {\tt{mil3}}, {\tt{mil5}} need in general less functional evaluations, but require the availability of a primitive $F$.

In Table \ref{table1} we show a couple of experiments in order to test the code for automatic computation of the alternating series described in Appendix \ref{appendix A1}. 
Given a certain tolerance, the idea is to use formula (\ref{stima errore}) to select $M$ with $C$ as in (\ref{costante C}).

In these tests,
\begin{itemize}
\item $f_4(x)=\frac{1}{(x+1)^{1/2}+1}$, and the reference value is $S^{(4)}=0.27828255712579120$;
\item $f_5(x)=\frac{1}{\log{(x+1)}+1}$, and the reference value is $S^{(5)}=0.66566699137697505$.
\end{itemize}
We report for the tolerances $10^{-4}$, $10^{-8}$, $10^{-12}$ the number of functional evaluations $m$ and the actual absolute error.
In both the examples estimates are accurate, with an error usually a little smaller than the required tolerance.

\begin{table}[ht]
\begin{center}
\begin{tabular}{cccccc}
\toprule
&  \multicolumn{2}{c}{$f(x) = \frac{1}{(x+1)^{1/2}+1}$} &  \multicolumn{2}{c}{$f(x) = \frac{1}{\ln(x+1)+1}$} \\
tolerance &  $m$ & error & $m$ & error  \\
\midrule
$1e-4$  & $6$ & $5.4e-5$ & $7$ & $6.8e-5$  \\
$1e-8$  & $18$ & $3.0e-9$ & $20$ & $5.3e-9$  \\
$1e-12$  & $38$ & $3.1e-13$ & $41$ & $6.1e-13$  \\
\bottomrule
\end{tabular}
\end{center}
\caption{Results of automatic summation with prescribed tolerance, in which $m$ is the total number of function evaluations and the error is computed by means of a reference solution.}\label{table1}
\end{table}

\section{Conclusions}

In this work we have presented a very simple algorithm based on the sinc rule for the summation of alternating series that can be written by means of the Abel formula. An accurate a-priori error estimation is derived and used to design the code for automatic evaluation of the sum of the series with prescribed accuracy. Although we absolutely do not claim that this method is the fastest with respect to the number of function evaluations (as reported in Figures \ref{fig1}-\ref{fig2}), we emphasize that it does not require additional information on the function $f$. Moreover, unlike standard situations, the a-priori error estimate (\ref{stima errore}) does not involve further analysis on the function $f$, such as derivatives or location of singularities.
  
\section*{Acknowledgements} 
Work partially supported by the DOR funds of the University of Padova and
by the INdAM-GNCS 2026 projects. This research has been accomplished within
the Community of Practice “Green Computing” of the Arqus European University
Alliance, the RITA “Research ITalian network on Approximation”.

\appendix

\section{Matlab codes}

In this appendix we describe the Matlab codes that we have used to perform our numerical experiments. They can be downloaded as open-source software at {\cite{GITDNS}}.

\subsection{The routine {\tt{sinc\_summation}}} \label{appendix A1}

The routine {\tt{sinc\_summation}}, given the function $f$ defining the alternate series $A(f)$, and a tolerance {\tt{tol}}, determines {\tt{s}} via the simple algorithm described in this paper, employing the sinc rule, i.e.,
\begin{equation}\label{437}
{\mbox{\tt{s}}}=f(0)-\frac{1}{2}f(1) +h \left( \frac{\Im(f(i\delta+1)}{2 \delta \pi}+\sum_{k=1}^M s_{\tilde{f}}(kh) \right),
\end{equation}
so that
$$
|{\tt{s}} - A(f)| \leq {\tt{tol}}.
$$
The value of $M$ is automatically selected by using (\ref{stima errore}) with the constant $C$ as in (\ref{costante C}).
The code also provides the number of functional evaluation required by such a summation algorithm, that is, $M+3$. 

Alternatively, given $M$, it computes {\mbox{\tt{s}}} as in (\ref{437}), as well as the error estimate (\ref{stima errore}).

\subsection{The routine {\tt{milovanovic\_summation}}}
Concerning the summation algorithm by G.V. Milovanovic, in \cite{M94} it is proved that if 
$$
A_m(f)=\sum_{k=m}^{\infty} (-1)^k f(k)
$$
and $F$ is a primitive of $f$ such that
\begin{itemize}
\item $F$ is holomorphic in the region $\{z \in {\mathbb{C}}: Re(z) \geq \alpha, m-1 < \alpha < m \}$,
\item $\lim_{|t| \rightarrow \infty} e^{-c|t|} f(x+i t/\pi)=0$, uniformly for $x \geq \alpha$,
\item $\lim_{x \rightarrow \infty} \int_{-\infty}^{+\infty} e^{-|t|} f(x+i t/\pi) dt=0$,
\end{itemize}
then,
$$
A_m(f)=\int_{-\infty}^{+\infty} {\boldsymbol{\psi}}(\alpha,t/(2\pi)) \sinh(t/2) \frac{e^{-t}}{(1+e^{-t})^2} dt,
$$
where
$$
{\boldsymbol{\psi}}(x,y) = \frac{(-1)^m}{2i} [F(x+iy)-F(x-iy)].
$$
Next, if $\{\lambda_{\eta}\}_{\eta=1\ldots,n}$, $\{\tau_{\eta}\}_{\eta=1\ldots,n}$ are the weights and nodes of the $n$-point Gaussian rule w.r.t. the weight function
$$
w(x)=\frac{\sinh(x)}{\cosh^2(x)}, \quad x \in (-\infty,\infty),
$$
then,
\begin{align*}
A_m(f)&=\sum_{\eta=1}^n \lambda_{\eta} {\boldsymbol{\psi}} \left(m-\frac{1}{2},\frac{1}{2\pi} \tau_{\eta} \right) \sinh\left(\frac{\tau_{\eta}}{2\pi} \right) + R_n({\boldsymbol{\psi}}) \\    
&=A_{m,n}(f)+\mathcal{R}_n({\boldsymbol{\psi}}),
\end{align*}
where $\mathcal{R}_n$ is the n-point Gaussian rule remainder.

Once $A_{m,n}(f)$ is available, we immediately get 
$$
A(f) \approx \sum_{k=0}^{m-1} (-1)^k f(k) + A_{m,n}(f).
$$
To this purpose, we have implemented in Matlab a procedure for the computation of a $n$-point Gaussian rule w.r.t. the weight function $w$, via Golub-Welsch algorithm. 

In particular, taking into consideration that $w$ is numerically null in $[350,+\infty)$, the recurrence coefficients are determined by the routine {\tt{r\_milovanovic}}, via the application of a discretization procedure (Stieltjes algorithm, see e.g. \cite[p.95]{GAUBOOK2004}) based on Gauss-Legendre quadrature rule in the intervals $I_1=[0,21.875]$ and $I_2=[21.875,350]$ with algebraic degree of exactness respectively $n_1+(2n-1)$ and $n_2+(2n-1)$, where $n_1=104$, $n_2=77$. 
The quantities $n_k$, $k=1,2$ depend on the fact that the weight function $w$ can be approximated at machine precision, in the interval $I_k$, $k=1,2$, by polynomials of such degrees, as it can be easily checked using Chebfun environment (see \cite{chebfun}). 
Next, the factor $2n-1$ comes from the requirement that a certain scalar product of two general polynomials of degree at most $n$ and $n-1$ must be computed exactly. 

%\textcolor{black}{ELEONORA 21.08: non siamo esperti: serve dire cos'è ADE?}

In the routine {\tt{milovanovic\_summation}}, using the tools just described, we first compute a $n$-points Gaussian rule w.r.t. $w(x)=\frac{\sinh(x)}{\cosh^2(x)}$, $x \in {\mathbb{R}}$ and then apply the summation technique suggested above.

\subsection{The routines {\tt{demo\_summation}}, {\tt{demo\_summation\_tables}}}

All the figures in this work are achieved using the routine {\tt{demo\_summation}}, that defines a test from a battery of examples, perform the summation methods increasing the number of functional evaluations and finally plots the absolute errors w.r.t. a reference result. For the {{sinc}} summation algorithm it also displays the error estimate (\ref{stima errore}).

\textcolor{black}{Concerning the Table \ref{table1}, the numerical results are obtained by calling the function {\tt{demo\_summation\_tables}}, in which we apply {\tt{sinc\_summation}}, varying the tolerance.}

\printbibliography

\end{document}